\documentstyle[12pt,twoside,amssymb]{article}
\def\Rightheadtext{Divergences as a Grading of the Formal
Variational Calculus} 
\def\Leftheadtext{Vladimir O. Soloviev} 
\def\Keywords{hamiltonian formalism, field theory, Poisson brackets,
boundary terms} 
\def\Classification{1991 {\it  Mathematics Subject
Classification}: Primary 58F05; Secondary 70G50, 58G20} 
\def\inprod{\mathop{\kern -0.05em\raise -0.1em\hbox{%
  \vrule height 0.03em width 0.6em depth 0em%
  \vrule height 0.7em width 0.03em depth 0em}\kern 0.1em}\nolimits}
\def\d{\mbox{\sf d}}
\def\Rn{{\Bbb R}^n}
\def\Tr{\mathop{\rm Tr}\nolimits}
\def\keywords#1{\footnotetext{{\it Key words and phrases:\/} #1}}
\def\classification#1{\footnotetext{#1}}
\def\Abstract#1{{\leftskip1cm\rightskip1cm\footnotesize\noindent
{\sc Abstract.} #1 \par}}
\font\bbf=cmbxti10 scaled\magstep1
\makeatletter
\def \section{\removelastskip\@startsection{section}{1}
  {\z@}
  {-3.5ex plus-1ex minus-.2ex}
  {2.3ex plus.2ex}
  {\def\@svsec{\thesection.\space}\noindent\hfil\sc}}
\def \subsection{\removelastskip\@startsection{subsection}{1}
  {\z@}
  {-3.5ex plus-1ex minus-.2ex}
  {2.3ex plus.2ex}
  {\def\@svsec{\thesubsection.\space}\noindent\hfil\sc}}
\makeatother

\def \thebibliography#1{\section*{References}
  \frenchspacing
  \raggedbottom
  \footnotesize
  \list{[\arabic{enumi}]}{\settowidth\labelwidth{[#1]}
    \leftmargin\labelwidth \advance\leftmargin\labelsep
    \itemsep 0pt \parsep 0pt
    \usecounter{enumi}}
  \def \newblock{\hskip .1em plus .3em minus .07em}
  \sloppy\clubpenalty10000\widowpenalty10000}

\begin{document}
\setcounter{page}{1}
\def\Refs{\list{}{\topsep-20pt\itemsep0pt\parsep0pt}\item[]}
\let\endRefs=\endlist
\def\References{\begin{center}
{\sc References}
\end{center}
\vspace*{-17.5pt}
\footnotesize}
\newcount\firstpage
\firstpage\thepage
\makeatletter
\renewcommand{\@oddhead}
   {\ifnum\thepage=\firstpage{}\else{\hfil\sc
    \lowercase\expandafter{\Rightheadtext} \hfil\thepage}\fi}
\renewcommand{\@evenhead}
   {\ifnum\thepage=\firstpage{}\else{\thepage\hfil\sc
\lowercase\expandafter{\Leftheadtext} \hfil}\fi}
\long\def\@makefntext#1{\parindent 1em\noindent
                   \ifnum\@thefnmark=0 \hbox to1em{\hss}
                   \else\hbox to1.8em{\hss$\m@th^{\@thefnmark}$}\fi #1\message{footmark=\@thefnmark}}
\makeatother

\begin{flushright}
\bf International Conference on\\
\bbf Secondary calculus and cohomological Physics,\\
\bf Moscow, August 1997
\end{flushright}

\classification{\Classification}


\vskip.5cm

\begin{center}
\uppercase{\bf Divergences as a Grading of the Formal Variational
Calculus}\\ 
\vspace{1.084cm}
\sc Vladimir O. Soloviev\\  
\vspace{.373cm}
\rm Institute for High Energy Physics,\\
142284 Protvino, Moscow Region, Russia\\
{\it E-mail:\/} vosoloviev@mx.ihep.su \\ 
\end{center}
\vskip.576cm

\Abstract{It is shown that the new formula for the 
field theory Poisson brackets  arise naturally in the
extension of the formal variational calculus incorporating divergences. The
linear spaces of local functionals, evolutionary vector fields, functional
forms, multi-vectors and differential operators become graded with respect to
divergences. The bilinear operations, such as the action of vector fields on
functionals, the commutator of vector fields, the interior product of  forms
and vectors and the Schouten-Nijenhuis bracket are compatible with the
grading.  A definition of the adjoint graded operator is proposed and
antisymmetric operators are constructed with the help of boundary terms.
The fulfilment of the Jacobi identity for the new Poisson brackets is shown
to be equivalent to vanishing of the Schouten-Nijenhuis bracket of the
Poisson bivector with itself. It is
demonstrated, as an example, 
that the second structure of the Korteweg-de Vries equation is
not Hamiltonian with respect to the new brackets until special
boundary conditions are prescribed.}

\section {Introduction}
The Hamiltonian formalism of classical mechanics may serve as the
ideal model illustrating the harmony of physics and mathematics.
During the last 20 years it was realized that a number of its mathematical
constructions, for example, the Schouten-Nijenhuis bracket\cite{Nij},
could be extrapolated to field theory\cite{Olv,Dorf}. This
made the search for new nonlinear integrable models much easier.
Even more general constructions uniting the Schouten-Nijenhuis bracket
with the Frolicher-Nijenhuis bracket were 
considered by A.~Vinogradov\cite{Vin}.

But the  methods mentioned above (and usually called 
the formal variational calculus\cite{GD}) have some restrictions 
arising from  boundary conditions which
should allow free integration by parts. As a rule, the periodic boundary
conditions or the rapid decay of fields at spatial infinity are necessary.  
Of course, these are not all the  physically interesting cases. For example,
the Coulomb potential in electrodynamics does not tend to zero rapidly
enough. The similar behaviour is typical for Yang-Mills and gravitation
fields. Nontrivial boundary problems arise also 
in the motion of material continua.

We consider the Hamiltonian
treatment of  nontrivial boundary problems
as an interesting field of research where there is some place for
new approaches and results.
The field theory
Poisson brackets which fulfil the Jacobi identity
under arbitrary boundary conditions have been proposed in \cite{Sol1}.  
Here we
extend the formal variational calculus to the most general case when no one
boundary term arising in integration by parts can be discarded.  We hope to
present some physical applications of the methods developed here in the future.
Let us also say that the interest to the role of divergences in field theory
is vivid now as can be seen from related papers \cite{BarnHenn, BarnStash,
Dick}.

As an example, illustrating the nonstandard nature of the problems to be
considered, let us remind the history of the long discussion on the role of
surface integrals in the canonical formulation of General Relativity.  During
about 15 years Arnowitt, Deser and Misner\cite{ADM}, Dirac\cite{Dirac},
Higgs\cite{Higgs}, Schwinger\cite{Schwinger}, DeWitt\cite{DeWitt}, Regge and
Teitelboim\cite{RT} were involved in it. The solution obtained in the last
work\cite{RT}  serves as the paradigm for the treatment of similar
problems up to now. It has been proposed to work with the special class of the
so-called differentiable or admissible functionals.  These
functionals are defined by the requirement that their variation should not
have surface contributions under the prescribed boundary conditions.  The
Poisson brackets for these functionals are the standard ones, i.e., they are
just the same as given in the formal variational calculus
\[ 
\{F,G\}=\int\limits_{\Omega}\left({{\delta F}\over{\delta
q^A(x)}} {{\delta G}\over{\delta p_A(x)}}-{{\delta G}\over{\delta q^A(x)}}
{{\delta F}\over{\delta p_A(x)}}\right)d^nx,
\]
but now nonzero surface contributions are allowed.

Here a natural question to ask is: do these brackets fulfil the standard
axiomatic requirements, i.e., the Jacobi identity and the existence of the
Poisson algebra on this space of admissible functionals? For the infinite
domain of integration and the asymptotic boundary conditions the affirmative
result for the second requirement was obtained by Brown and
Henneaux\cite{BH}. The first requirement was partially analysed by us and in
the case treated above the answer is also positive.

It is more difficult to study the finite domain case. Let us take as a second
example the motion of a fluid or plasma. It was shown by Lewis,
Marsden, Montgomery and Ratiu\cite{LMMR} that the Jacobi identity for
the standard Poisson bracket can be violated even in the case of fixed
boundary, and so the Poisson brackets should be modified by surface
terms. In the free boundary case it turns out to be natural to extend the
space of admissible  functionals so that their variation
could include nonzero surface contributions. But according to\cite{LMMR}
the presence of nonzero term with $\delta q^A$ in the boundary integral
requires the absence of the corresponding term with $\delta p_A$
and vice versa. A new formula for Poisson brackets arises as a result of a
generalization of the variational derivative which is now allowed to contain a
surface contribution 
\[
\delta H=
\int\limits_{{\Omega}}{{\delta^{\wedge}H}\over{\delta q^A}}
\delta q^A d^nx +\oint\limits_{\partial {\Omega}}
 {{\delta^{\vee} H}\over {\delta q^A}}
\delta q^A \vert _{\partial {\Omega}} dS+
\int\limits_{{\Omega}}{{\delta^{\wedge}H}\over{\delta p_A}}
\delta p_A d^nx +\oint\limits_{\partial {\Omega}}
 {{\delta^{\vee} H}\over {\delta p_A}}
\delta p_A \vert _{\partial {\Omega}} dS.
\]
Unfortunately, it is not quite clear whether the
Poisson bracket of the two functionals, admissible in this new
sense, will be admissible functional itself.

As a third example, we would like to attract reader's attention
to consequences of the noncommutativity of the standard variational
derivatives, i.e., the Euler-Lagrange derivatives. This point was
discussed  formally  in publications by
I.~Anderson\cite{And76,And78}
and Aldersley\cite{Ald}. We faced with the problem independently, in the
course of studying  surface terms in the Poisson algebra of the Ashtekar
formalism of canonical gravity\cite{Sol92}. There it was found  that
transformations of the type
\[
q^A(x)\to q^A(x),\qquad p_A(x)\to
p_A(x)+{{\delta F[q]}\over{\delta q^A(x)}},
\]
were noncanonical ones 
if surface terms were not ignored.  Tracing the connection of
this calculation with the standard calculations with
$\delta$-functions\cite{Ash}, we have found that the
correspondence could be restored by introduction of $\theta_{\Omega}$ ---
the characteristic function of the domain $\Omega$
\[
\theta_{\Omega}(x)= \cases {
1 &  if $x\in\Omega$; \cr
0  & otherwise. \cr }
\]
Then the standard relations
\[
\left({{\partial}\over{\partial x^i}}+{{\partial}\over{\partial y^i}}\right)
\delta(x,y)=0,
\]
should be modified as
\[
\left(\theta_{\Omega}(x){{\partial}\over{\partial
x^i}}+\theta_{\Omega}(y){{\partial}\over{\partial y^i}}\right)\delta(x,y) =
-{{\partial\theta_{\Omega}}\over{\partial x^i}}\delta(x,y),
\]
where the usually discarded surface term is preserved.

All the above examples lead us to the necessity to extend  the formal
variational calculus onto total divergences. This extension consists in the
introduction of a new grading for the linear spaces of local functionals,
vector fields, functional forms, multi-vectors and differential operators.
To come back to the standard case one should put $\theta_{\Omega}(x)\equiv 1$
in $\Rn$.

The extension of the formal variational calculus naturally incorporates
the new definition of local functionals (not modulo divergences) and of
their differentials (as a full variation, not fixed on the boundary).
The Poisson bivectors are now defined in a more general way and they can
include boundary contributions. The definition of pairing (interior product)
is also revised and now the trace of two differential operators is used
for it, so the
pairing becomes compatible with the grading. The Poisson bracket found more
or less heuristically in \cite{Sol1} arises now on the base of the geometrical
constructions as
\[
\{F,G\}=\d G\inprod\d F\inprod\Psi,
\]
where $\Psi$ is the Poisson bivector.

   We show here that the Jacobi identity for the new Poisson brackets can
be verified without the long calculations of binomial sums used in \cite{Sol1}.
 Its
fulfilment is equivalent to the vanishing of the Schouten--Nijenhuis bracket
of the Poisson bivector with itself. And in its turn this condition can be
easily tested along with the procedure proposed by  Olver\cite{Olv} with a
minimal modification of it. More attention than in \cite{Sol1} is paid here to
nonultralocal Hamiltonian operators with nonconstant coefficients,
because now the technical obstacles are removed. It turns out to be
that not all operators which are Hamiltonian with respect to the standard
brackets remain Hamiltonian in relation to the new brackets. 
For example, the second structure of the Korteweg--de Vries equation
is not automatically Hamiltonian with respect to the new formalism.
In this respect it strongly differs from the first KdV structure.

The content of this work is as follows.
In Section 2 we introduce the grading for the
local functionals, and the evolutionary vector fields. In Section 3 the
differential, the functional $m$-forms, the interior product of vector
fields onto $m$-forms and Lie derivative are defined.
Section 4 deals with graded differential operators and their adjoints. In
Section 5 we discuss multi-vectors and the  Schouten-Nijenhuis bracket. It
is shown how 1-vectors and evolutionary vector field spaces are isomorphic.
Section 6 contains the general definition of the Poisson bracket, its
realization in this formalism, the definition of Hamiltonian
vector fields and the theorem on connection between the commutator of two
Hamiltonian vector fields and the 
Poisson bracket of corresponding Hamiltonians.
All constructions are illustrated by an example: the first Hamiltonian
structure of the Korteweg-de Vries equation.
The proof of Jacobi
identity is postponed until Section 7. This new proof is applicable for
all local Poisson brackets and so supersedes the proofs given earlier in 
\cite{Sol1}. In
the same time it is easy to compare this proofs because notations are the same.
At last, in Section 8 we consider two
examples of the non-ultralocal Poisson brackets with field dependent
coefficients (this class of brackets we were unable to study by the
methods of \cite{Sol1}). We show that  the second structure of Korteweg-de Vries
equation is not Hamiltonian if surface terms are not ignored, whereas  the
2-dimensional flow of the ideal fluid is described by Hamiltonian structure.
This points to nontrivial character of the generalization of the formal
variational calculus suggested here.
A short summary is given in Conclusion and the points remaining unclear are
mentioned.  We hope to continue this research by 
the detailed treatment of the boundary conditions role and 
applications to some physical problems.

As a rule,  we use below the same notations as in \cite{Sol1} except a
change of the notation for the Fr\'echet derivative  from $D_{f}$ to $f'$
and omitting the sign of summation according to the Einstein rule.
We find convenient to represent integrals over finite domain $\Omega$
through integrals over the infinite space $\Rn$ by inserting
into all integrands  the characteristic function
$\theta_{\Omega}$. Then the formalism seems closer to the standard
formal variational calculus where local functionals and functional forms
are defined modulo divergences. But the formal divergences that we
discard here are integrated to zero under arbitrary
conditions on the boundary of the finite domain, whereas real
divergences are incorporated into graded structures.
All the operations introduced below are compatible both with
discarding formal total divergences
(if one object is a formal divergence than the result of operation
is also formal divergence) and with the grading
(i.e., the same is valid for real divergences).
Extension of the space of differential operators by admitting their grading
permits
to use the concept of adjoint operator. So, antisymmetric operators can now
be constructed and the Poisson bracket formulas become more compact, than in
\cite{Sol1} though their content is the same. Nevertheless, in the proof of 
the Jacobi
identity we prefer to use the old notations to make easier the  comparison with
the not so general  proofs of \cite{Sol1}.

Henceforth we consider the space $\Rn$
and the multi-index notations $J=(j_1,...,j_n)$ where
$j_i\geq 0$
\[
\phi_A^{(J)}=
{{\partial^{|J|}\phi_A}\over{\partial^{j_1}x^1...{\partial^{j_n}x^n}}},
\qquad |J|=j_1+...+j_n.
\]
The Fr\'echet derivative is defined as
\begin{equation}
f'_A=\frac{\partial f}{\partial\phi_A^{(J)}}D_J,\label{eq:frechet}
\end{equation}
where
\[
D_i={{\partial}\over{\partial x^i}}+\phi_A^{(J+i)}{{\partial}\over
{\partial\phi_A^{(J)}}}, \quad D_J=D_1^{j_1}...D_n^{j_n},\quad D_i^0=1,
\quad D_i^{-1}=0.
\]
Binomial coefficients for multi-indices are
\[
{J \choose K}={j_1\choose k_1}\cdots{j_n\choose k_n},
\]
\[
{j \choose k}= \cases {
j!/(k!(j-k)!) & if  $ 0\le k \le j$; \cr
0  & otherwise. \cr }
\]
With the help of them we introduce
the so-called higher Eulerian operators \cite{Olv,Ald,KMGZ}
\begin{equation}
E^J_A(f)=(-1)^{|K|+|J|}{K\choose
J}D_{K-J}{{\partial f} \over{\partial\phi_A^{(K)}}}.\label{eq:higher}
\end{equation}

\section{Local functionals and evolutionary vector fields}
Let us start with notions from the theory of graded spaces as they are
given in Ref.~\cite{Dorf}.
A {\it  grading} in linear space $L$ is a decomposition of it into direct sum
of subspaces, with a special value of some function $p$ (grading function)
assigned to all the elements of any subspace.

Below the function $p$ takes its values in the set of all positive
multi-indices $J=(j_1,\dots,j_n)$ and so,
\[
L=\bigoplus\limits_{J=0}^{\infty} L^{\langle J\rangle}.
\]
Elements of each subspace are called homogeneous.

A bilinear operation $x,y\mapsto x\circ y$, defined on $L$, is said to be
{\it compatible with the grading} if the product of any homogeneous elements
is also homogeneous, and if
\[
p(x\circ y)=p(x)+p(y).
\]
Now let us turn to the concrete structures.

There are two ways to write a local functional: as the integral
of a smooth function $f^{\langle 0\rangle}\left(\phi^{(K)}_A(x)\right)$ of
fields and their derivatives up to some finite order over the
prescribed domain $\Omega$ in $\Rn$, or as the integral over all the space
$\Rn$ but with the characteristic function of the domain $\theta_{\Omega}$
included into the integrand
\begin{equation}
F=\int\limits_{\Omega}f^{\langle 0\rangle}\left(\phi^{(K)}_A(x)\right)d^nx
\equiv
\int\limits_{\Rn}
\theta_{\Omega}f^{\langle 0\rangle}d^nx. \label{eq:locfunc}
\end{equation}
As in \cite{Sol1}, let us denote the space of local functionals as $\cal A$.
Here we shall call the expression given above
the {\it canonical form of a local functional}. We formally extend that
definition by allowing local functionals to be written as follows
\[
F=\int\limits_{\Rn} D_J\theta_{\Omega}f^{\langle J\rangle}
\left(\phi^{(K)}_A(x)\right)d^nx\equiv \int\theta^{(J)}f^{\langle J\rangle}d^nx
\equiv\int fd^nx,
\]
where only a finite number of terms is allowed. Here and below we simplify
the notation for derivatives of $\theta$ and remove $\Omega$.
All the integrals without the domain of integration shown explicitly are
integrals over $\Rn$, below we shall omit $d^nx$. Of course, any functional can
be transformed to the above form (\ref{eq:locfunc}), exclusively used in \cite{Sol1},
through integration by parts
\[
F=\int\theta\tilde f^{\langle 0\rangle}\equiv
\int\limits_{\Omega}\tilde f^{\langle 0\rangle},
\]
where
\[
\tilde f^{\langle 0\rangle}= (-1)^{|J|}D_Jf^{\langle J\rangle}.
\]
Evidently, the formal integration by parts over infinite
space $\Rn$ changes the grading.
It will be clear below that the general
situation is the following: from one side we have the compatibility of all
the bilinear operations with the grading and from the other side --- with
the formal integration by parts.  So, basic objects (local functionals etc.) are
defined as equivalence classes modulo formal divergences (i.e., divergences
of expressions containing $\theta$-factors) and the unique decomposition into
the homogeneous subspaces with the fixed grading function can be made only
for representatives of these classes.

We call  expressions of the form
\[
\Psi= \int\theta^{(J)}
D_K\psi^{\langle J\rangle}_A\frac{\partial}{\partial\phi_A^{(K)}}\equiv
\int\theta^{(J)}\psi^{\langle J\rangle}\equiv\int\psi,
\]
the {\it evolutionary vector fields}.
The action of the evolutionary vector field on a local functional is
given by the expression
\begin{equation}
\Psi F= \int\theta^{(I+J)}
D_K\psi^{\langle J\rangle}_A\frac{\partial f^{\langle I
\rangle}}{\partial\phi_A^{(K)}}\equiv\int\theta^{(I+J)}
\psi^{\langle J\rangle}f^{\langle I\rangle}\equiv\int\psi f.\label{eq:mult}
\end{equation}
It is a straightforward calculation to check that this operation is
compatible with the formal integration by parts, i.e.
\[
\psi Df= D(\psi f),
\]
as it is in the standard formal variational calculus. This relation is,
of course, valid for integrands.

It is easy to see that the evolutionary vector field with coefficients
\[
\psi^{\langle J\rangle}_A=
 D_L\xi_B^{\langle I\rangle} \frac{\partial\lambda_A^{\langle J-I\rangle
}}{\partial \phi_B^{(L)}}-
D_L\lambda_B^{\langle I\rangle} \frac{\partial\xi_A^{\langle J-I\rangle}}
{\partial \phi_B^{(L)}}
\]
can be considered as the {\it  commutator of the evolutionary vector fields}
 $\Xi$ and $\Lambda$
\[
\Psi F=[\Xi,\Lambda]F=\int\biggl(\xi(\lambda f)-\lambda(\xi f)\biggr),
\]
with the Jacobi identity fulfilled for the commutator operation,
and so these vector fields form a Lie algebra.

Let us comment upon  the representation of the evolutionary vector fields
as integrals, which is different from the traditional notations.

The formal variational calculus\cite{GD} operates with the local functionals
which are expressed by single integrals of functions of the specified
class, for example, infinitely differentiable ones. The functional forms
and multi-vectors are expressed by similar integrals. The pairing of two
such objects gives us a single integral again.

At the same time, another notations are widespread, especially in physical
literature, which use $\delta$-function and its derivatives. Then a
result of the pairing of two single integrals is understood as a double
integral. But as this double integral contains the $\delta$-function,
it always can be converted into the single one.

This convertation of a double integral into the single one with the help
of $\delta$-function is trivial when no boundary terms could arise.
The subject of this work is just a study of the opposite case. The new
rule is necessary here and it have been proposed in \cite{Sol1} as Rule~4.2
\begin{equation}
\int\limits_{\Omega}
\int\limits_{\Omega} f(x)g(y)D_J^{(x)}D_K^{(y)}\delta (x,y)=
\int\limits_{\Omega} D_KfD_Jg.\label{eq:rule42}
\end{equation}
In
this article we give really a new and equivalent form of these rules which
help to avoid the usage of double integrals at all.

The concept of vector field appeared initially in the course of studying
the evolutionary differential equations and their symmetries. In the formal
variational calculus\cite{GD} functionals are, in fact, replaced by
equivalence classes of functions, and so the action of evolutionary
vector fields onto local functionals is replaced by their action on functions
\[
\psi f=
D_K\psi_A\frac{\partial f}{\partial\phi_A^{(K)}}.
\]
However, to represent functionals by integrals
and to require that the result of the action of an evolutionary vector
field onto a local functional is a local functional, i.e.
an integral, it is absolutely natural  to represent the evolutionary
vector fields also as integrals
\[
\Psi = \int\limits_{\Omega}
D_K\psi_A(x)\frac{\partial }{\partial\phi_A^{(K)}(x)}\equiv\int\psi,
\]
in combination with the standard rule
\begin{equation}
\frac{\partial\phi^{(J)}_A(y)}{\partial\phi^{(K)}_B(x)}=
\delta (x,y)\delta_{AB}\delta_{JK}.\label{eq:rule}
\end{equation}

Another argument supporting our notations is the equivalence between
evolutionary vector fields and 1-vectors, which is demonstrated for
the standard formal variational calculus in book\cite{Olv} and also
for the graded case in Section 5 of this article.  1-vectors as a partial
case of multi-vectors are always written as integrals.

Apart from the notational revision we would like to mention a
new feature in our treatment of the vector fields:
now they are not derivations when applied to standard functions, but only
to the graded ones.
Of course, in the traditional approach the vector fields are
not derivations when applied to functionals, because their
multiplication is not defined. But  these vector
fields, traditionally written without the integral sign,
are derivations of functions. This property is
partially lost here. It can be restored formally if we  
consider integrands containing $\theta$ as functions and take a relation
\begin{equation}
D_I\theta\times D_J\theta =D_{I+J}\theta.\label{eq:rule54}
\end{equation}
as a definition of their multiplication.

In this context, formula (\ref{eq:mult}), introduced as a
definition,  can be interpreted  also as a consequence of the standard
relation (\ref{eq:rule}) and a new definition (\ref{eq:rule54}).

Therefore, it is evident that our ``rule for multiplication of
distributions'' taken from \cite{Sol1}, i.e., 
equation (\ref{eq:rule54}) is nothing more than
another way to define the pairing compatible with the introduced grading.

At last, let us mention the possibility to use other notations in this
formalism. It is, of course, possible to avoid $\theta$-functions and to
use integrals over the domain $\Omega$ only. Then any local
functional can be given as
\[
F=\int\limits_{\Omega}D_Jf^{\langle\langle J\rangle\rangle},
\]
where
\[
f^{\langle\langle J\rangle\rangle}=(-1)^{|J|}f^{\langle J\rangle},
\]
with analogous rewriting of the other objects. Correspondingly,
equation (\ref{eq:mult}) will be written as
\[
\Psi F= \int\limits_{\Omega}D_{I+J}\left(
D_K\psi^{\langle\langle J\rangle\rangle}_A\frac{\partial
f^{\langle\langle I
\rangle\rangle}}{\partial\phi_A^{(K)}}\right).
\]

\section{Differentials and functional forms}
The {\it differential of a local functional}
is simply the first variation of it
\[
\d F= \int\theta^{(J)}
\frac{\partial f^{\langle J\rangle}}{\partial\phi_A^{(K)}}\delta\phi_A^{(K)}
\equiv\int\theta^{(J)}\d f^{\langle J\rangle}\equiv\int\d f,
\]
here and below $\delta\phi_A^{(K)}=D_K\delta\phi_A $.
It can also be expressed through
the Fr\'echet derivative (\ref{eq:frechet})
or through the higher Eulerian operators (\ref{eq:higher})
\[
\d F= \int\theta^{(J)}{f^{\langle J\rangle}}'(\delta\phi)=
 \int\theta^{(J)}
D_K\left( E^K_A(f^{\langle J\rangle})\delta\phi_A\right) .
\]
This differential is a special example of functional 1-form.
A general functional 1-form can be written as
\[
\Sigma =  \int\theta^{(J)}\sigma^{\langle J\rangle}_{AK}
\delta\phi_A^{(K)}\equiv\int\theta^{(J)}\sigma^{\langle J\rangle}\equiv\int\sigma.
\]
Of course, the coefficients $\sigma^{\langle J\rangle}_{AK}$ are not unique
since we can make formal integration by parts. Let us call the following
expression the {\it canonical form of a functional 1-form}
\[
\Sigma= \int\theta^{(J)}\sigma^{\langle J\rangle}_A\delta\phi_A.
\]
Analogously, we can define {\it functional} $m$-{\it forms}
as integrals or equivalence classes modulo formal
divergences of vertical $m$-forms
\[
\Sigma =\frac{1}{m!}
 \int\theta^{(J)}\sigma^{\langle J\rangle}_{A_1K_1,\dots,A_mK_m}\delta
\phi_{A_1}^{(K_1)}\wedge\dots\wedge\delta\phi_{A_m}^{(K_m)}
=
 \int\theta^{(J)}\sigma^{\langle J\rangle}=
 \int\sigma.
\]
Define the {\it pairing} (or the {\it interior product})
of an evolutionary vector field and 1-form as
\begin{equation}
\Sigma (\Xi)=\Xi\inprod\Sigma=
\int\theta^{(I+J)}\sigma^{\langle J\rangle}_{AK}D_K\xi_A^{\langle I\rangle}=
\int\theta^{(I+J)}\sigma^{\langle J\rangle}(\xi^{\langle I\rangle})=
\int\sigma(\xi).
\label{eq:pairing}
\end{equation}
The interior product of an evolutionary vector field and a functional $m$-form
will be given as follows
\[
\Xi\inprod\Sigma=
\frac{1}{m!} (-1)^{i+1}
\int\theta^{(I+J)}\sigma^{\langle J\rangle}_{A_1K_1,\dots,A_mK_m}
D_{K_i}\xi_{A_i}^{\langle I\rangle}\delta
\phi_{A_1}^{(K_1)}\wedge\dots
\]
\[
\dots\wedge\delta\phi_{A_{i-1}}^{(K_{i-1})}
\wedge\delta\phi_{A_{i+1}}^{(K_{i+1})}\wedge\dots\wedge
\delta\phi_{A_m}^{(K_m)}.
\]
Then a value of the $m$-form on the $m$ evolutionary vector fields
will be defined by the formula
\[
\Sigma (\Xi_1,\dots,\Xi_m)=\Xi_m\inprod\dots\Xi_1\inprod\Sigma.
\]
It can be  checked by straightforward calculation that
\[
(D\sigma)(\xi_1,\dots,\xi_m)=
D\left(\sigma(\xi_1,\dots,\xi_m)\right).
\]
The {\it differential of the} $m$-{\it form} which is given as
\[
\d\Sigma =\frac{1}{m!}
\int\theta^{(J)}\frac{\partial\sigma^{\langle J
\rangle}_{A_1K_1,\dots,A_mK_m}}
{\partial\phi_A^{(K)}}\delta\phi_A^{(K)}\wedge\delta
\phi_{A_1}^{(K_1)}\wedge\dots
\wedge\delta\phi_{A_m}^{(K_m)}=
\int\theta^{(J)}\d\sigma^{\langle J\rangle}=\int\d\sigma,
\]
satisfies standard properties
\[
{\d}^2=0
\]
and
\[
\d\Sigma(\Xi_1,\dots,\Xi_{m+1})=
\sum_i(-1)^{i+1}\Xi_i\Sigma(\Xi_1,\dots,
\hat\Xi_i,\dots,\Xi_{m+1})+
\]
\[
+\sum\limits_{i<j}(-1)^{i+j}\Sigma([\Xi_i,\Xi_j],\Xi_1,\dots,\hat\Xi_i,\dots,
\hat\Xi_j,\dots,\Xi_{m+1}).
\]
The {\it Lie derivative} of a functional form $\Sigma$
along the evolutionary vector field $\Xi$ can be introduced by the standard
formula
\[
L_{\Xi}\Sigma=\Xi\inprod\d\Sigma+\d \biggl(\Xi\inprod\Sigma\biggr).
\]

\section{Graded differential operators and their adjoints}
We call linear matrix differential operators of the form
\[
\hat I_{AB}=\theta^{(J)}
I^{\langle J\rangle N}_{AB}D_N
\]
the {\it graded differential operators}.

Let us call the linear differential operator $\hat I^{\ast}$ the {\it
adjoint} to $\hat I$ if for an
arbitrary set of smooth functions $f_A$, $g_A$
\[
\int f_A\hat I_{AB}g_B= \int g_A\hat
I^{\ast}_{AB}f_B.
\] 
For coefficients of the adjoint operator we can derive
the expression
\begin{equation}
I^{\ast\langle J\rangle
M}_{AB}=
(-1)^{|K|}{K\choose L}{K-L\choose M}
D_{K-L-M}I^{\langle J-L\rangle K}_{BA}.\label{eq:adj}
\end{equation}
It is easy to check that the relation
\[
\hat I_{AB}(x)\delta(x,y)=\hat I^{\ast}_{BA}(y)\delta(x,y)
\]
follows from Rule 4.2 of \cite{Sol1}. For example, we have
\begin{equation}
\left(\theta(x)\frac{\partial}{\partial
x^i}+\theta(y)\frac{\partial}{\partial y^i}\right) \delta
(x,y)=-\theta^{(i)}\delta (x,y).\label{eq:delta}
\end{equation}
In one of our
previous publications \cite{Sol92} we  tried to connect the appearance of
surface terms in Poisson brackets and the standard manipulations with the
$\delta$-function. The ansatz used there for
the above simplest example coincided with (\ref{eq:delta}) up to the sign.
The reason for this
difference laid in the other choice made there instead of Rule 4.2 of \cite{Sol1}.
That ansatz lead us to the standard Poisson brackets which are not
appropriate for nontrivial boundary problems.

Operators satisfying the relation
\[
\hat I^{\ast}=-\hat I
\]
will be called the {\it antisymmetric} ones. With the help of them it is possible
to express 2-forms (and also 2-vectors to be defined below) in the canonical
form
\[
\Sigma=\frac{1}{2}\int\delta\phi_A\wedge\hat I_{AB}
\delta\phi_B.
\]
It is clear that we can consider representations of functional forms
as decompositions over the basis derived as a tensor product
of $\delta\phi_A$, with the totally antisymmetric multilinear operators
\[
\hat\sigma= \theta^{(J)}\sigma^{\langle J\rangle}_{A_1K_1,\dots,A_mK_m}
\biggl( D_{K_1}\cdot,\dots,D_{K_m}\cdot\biggr)
\]
as coefficients of these decompositions.

\section{Multi-vectors, mixed tensors and Schouten-Nijenhuis bracket}
Let us introduce dual basis to $\vert\delta\phi_A\rangle$ by the relation
\begin{equation}
\left\langle\frac{\delta}{\delta\phi_B(y)},\delta\phi_A(x)\right\rangle
=\delta_{AB}\delta(x,y) \label{eq:dual}
\end{equation}
and construct by means of the tensor product a basis
\[
\frac{\delta}{\delta\phi_{B_1}(y)}\otimes\frac{\delta}{\delta\phi_{B_2}(y)}
\otimes\dots\otimes\frac{\delta}{\delta\phi_{B_m}(y)}.
\]
Then by using totally antisymmetric multilinear operators described in
the previous Section we can define the {\it functional $m$-vectors} (or {\it
multi-vectors})
\[
\Psi=\frac{1}{m!} \int\theta^{(J)} \psi^{\langle
J\rangle}_{B_1L_1,\dots,B_mL_m}D_{L_1}
\frac{\delta}{\delta\phi_{B_1}}\wedge\dots\wedge D_{L_m}
\frac{\delta}{\delta\phi_{B_m}}=
\int\theta^{(J)} \psi^{\langle
J\rangle}.
\]
Here a natural question on the relation between evolutionary
vector fields and 1-vectors arises.
Evidently, evolutionary vector fields lose their form when being
integrated by parts whereas 1-vectors preserve it.
Let us make a partial integration in the expression
of a general evolutionary vector field
\[
\Xi= \int\theta^{(J)}
D_K\xi^{\langle J\rangle}_A\frac{\partial}{\partial\phi_A^{(K)}}
\]
by removing $D_K$ from $\xi^{\langle J\rangle}_A$, then we get
\[
\Xi=
\int\xi^{\langle J\rangle}_A\theta^{(J+L)}
(-1)^{|K|}{K \choose L}D_{K-L}\frac{\partial}{\partial\phi_A^{(K)}}.
\]
It is easy to see that by using (\ref{eq:rule54}), i.e., Rule 5.4 from \cite{Sol1},
in the backward direction we can write
\[
\Xi= \int\left( \theta^{(J)}\xi_A^{\langle J\rangle}\right)\left(
\theta^{(L)}(-1)^{|L|}E^L_A\right)= \int\theta^{(J)}\xi^{\langle
J\rangle}_A
\frac{\delta}{\delta\phi_A},
\]
where the higher Eulerian operators (\ref{eq:higher}) and the full
variational derivative (Definition 5.1 of \cite{Sol1})
\[
\frac{\delta F}{\delta\phi_A}=\sum (-1)^{|J|}E^J_A(f)D_J\theta,
\]
are consequently used. Therefore, we arrive at the following Statement.

{\bf Statement 5.1}
{\it There is a one-to-one correspondence between the evolutionary vector
fields and the functional 1-vectors.
The  coefficients of 1-vector in the canonical form $\xi_A^{\langle J
\rangle}$ are equal to the characteristics
of the evolutionary vector field.}

It is not difficult to show that we can  deduce the pairing (interior
product) of 1-forms and 1-vectors and this pairing preserves the
identification.  Really, the definition of the dual basis (\ref{eq:dual}) and
(\ref{eq:rule54}), i.e., Rule 5.4 of \cite{Sol1}, permits us to derive that
\[
\Sigma(\Xi)=\Xi\inprod\Sigma= \int \int\theta^{(I)}(x)
\theta^{(J)}(y)\sigma^{\langle I\rangle}_{AK}(x)\xi^{\langle J\rangle
}_{BL}(y)
\left\langle D_L\frac{\delta}{\delta\phi_B(y)},D_K\delta\phi_A(x)\right
\rangle=
\]
\[
= \int\theta^{(I+J)}D_L\sigma^{\langle I\rangle}_{AK}
D_K\xi^{\langle J\rangle}_{AL}=
\int\theta^{(I+J)}\sigma^{\langle I\rangle}
(\xi^{\langle J\rangle})=\int\sigma(\xi)
= \int\theta^{(I+J)}{\Tr}(\sigma^{\langle I\rangle}
\xi^{\langle J\rangle}),
\]
and when 1-vector is in the canonical form (only $L=0$ term is nonzero)
this result coincides with Eq.(\ref{eq:pairing}).

This formula for the pairing will be exploited below also for interior
product of 1-vectors and $m$-forms or 1-forms and $m$-vectors.
Its importance comes from the fact that it is invariant under the formal
partial integration both in forms and in vectors, i.e., 
\[ 
({\rm D}\sigma)(\xi)={\rm D}(\sigma(\xi))=\sigma({\rm D}(\xi)).  \]
Evidently, it is the trace construction for convolution of differential
operators (as coefficients of tensor objects in the proposed basis)
that guarantees this invariance.

The interior product of 1-vector onto $m$-form and, analogously,
of 1-form onto $m$-vector is defined as
\begin{eqnarray}
\Xi\inprod\Sigma&=&\frac{1}{m!}  (-1)^{(i+1)}\int
\theta^{(I+J)}D_{K_i}\xi^{\langle I\rangle}_{A_iL}D_L
\Biggl(\sigma^{\langle J\rangle}_{A_1K_1,\dots,
A_mK_m}\delta\phi_{A_1}^{(K_1)}\wedge\dots\nonumber\\
&\dots&\wedge\delta\phi_{A_{i-1}}^{(K_{i-1})}\wedge
\delta\phi_{A_{i+1}}^{(K_{i+1})}\wedge\dots
\wedge\delta\phi_{A_m}^{(K_m)}\Biggr)=
(-1)^{(i+1)}\int
\theta^{(I+J)}\xi^{\langle I\rangle}\inprod\sigma^{\langle J\rangle}.
\end{eqnarray}
Then we also can define the value of $m$-form on $m$ 1-vectors (or,
analogously, $m$-vector on $m$ 1-forms)
\[
\Sigma (\Xi_1,\dots,\Xi_m)=
\Xi_m\inprod\dots\Xi_1\inprod\Sigma=
 \int\theta^{(J+I_1+\dots+I_m)}{\Tr}
\left( \sigma^{\langle J\rangle}
\xi_1^{\langle I_1\rangle}\cdots\xi_m^{\langle I_m\rangle}\right),
\]
where each entry of multilinear operator $\sigma$ acts only on the 
corresponding
$\xi$, whereas each derivation of the operator $\xi$ acts on the
product of $\sigma$ and all the rest of $\xi$'s.

It is possible to define the {\it differential of $m$-vector}
\[
\d\Psi =\frac{1}{m!} \int\theta^{(J)}
\frac{\partial\psi^{\langle J\rangle}_
{A_1K_1,\dots,A_mK_m}}{\partial\phi_B^{(L)}}\delta\phi_B^{(L)}D_{K_1}
\frac{\delta}{\delta\phi_{A_1}}\wedge\dots\wedge
D_{K_m}\frac{\delta}{\delta\phi_{A_m}},
\]
as an example of a mixed ${m \choose 1}$ object. Evidently, ${\d}^2\Psi=0$.

With the help of the previous constructions we can define the
{\it Schouten-Nijenhuis bracket}
\[
\bigl[ \Xi,\Psi\bigr]_{SN} =\d\Xi\inprod\Psi +
(-1)^{pq}\d\Psi\inprod\Xi
\]
for two multi-vectors of orders
$p$ and $q$. The result of this operation is $p+q-1$-vector and it is
analogous to the Schouten-Nijenhuis bracket in tensor analysis \cite{Nij}.
Its use in the formal variational calculus is described in
Refs.\cite{Olv,Dorf}.  However, in cited references this bracket is
usually defined for operators.  We can recommend Ref.\cite{Olv2} as an
interesting source for the treatment of the Schouten-Nijenhuis bracket of
multi-vectors.  Our construction of this bracket guarantees a compatibility
with the equivalence modulo divergences

\[ \bigl[ D\xi,\psi\bigr]_{SN} =D\bigl[ \xi,\psi\bigr]_{SN}= \bigl[ \xi,
D\psi\bigr]_{SN}.
\]
\medskip
{\bf Statement 5.2} {\it The
Schouten-Nijenhuis bracket of functional 1-vectors up to a sign coincides
with the commutator of the corresponding evolutionary vector fields.}
\medskip

{\it Proof.}
Let us take the two 1-vectors in canonical form without loss of generality
\[
\Xi= \int\theta^{(J)}\xi^{\langle J\rangle}_A
\frac{\delta}{\delta\phi_A},\qquad
\Psi= \int\theta^{(K)}\psi^{\langle K\rangle}_B\frac{\delta}{\delta\phi_B}
\]
and compute
\[
\bigl[ \Xi,\Psi\bigr]_{SN}=\d\Xi\inprod\Psi
-\d\Psi\inprod\Xi.
\]
We have
\[
\d\Xi= \int\theta^{(J)}{\xi^{\langle
J\rangle}_A}'(\delta\phi)\frac{\delta}{\delta
\phi_A}= \int\theta^{(J)}\frac{\partial\xi^{\langle J\rangle}_A}{\partial\phi^{(L)}
_C}\delta\phi_C^{(L)}\frac{\delta}{\delta\phi_A},
\]
and
\[
\d\Xi\inprod\Psi=- \int\theta^{(J+K)}\frac{\partial
\xi_A^{\langle J\rangle}}
{\partial\phi_B^{(L)}}D_L\psi_B^{\langle K\rangle}\frac{\delta}{\delta\phi_A}.
\]
Therefore, we obtain
\[
\bigl[ \Xi,\Psi\bigr]_{SN}=- \int\theta^{(J+K)}\left(
D_L\psi_B^{\langle K\rangle}\frac{\partial\xi_A^{\langle J\rangle}}
{\partial\phi_B^{(L)}}-
D_L\xi_B^{\langle K\rangle}\frac{\partial\psi_A^{\langle J\rangle}}
{\partial\phi_B^{(L)}}
\right)
\frac{\delta}{\delta\phi_A}=-[\Xi,\Psi],
\]
and the proof is completed.

\medskip
{\bf Statement 5.3} (Olver's Lemma \cite{Olv})
{\it
The Schouten-Nijenhuis bracket for two bivectors can be expressed
in the form}
\begin{equation}
\bigl[ \Lambda,\Psi\bigr]_{SN}=-\frac{1}{2} \int\
\xi\wedge\hat I'(\hat K\xi)\wedge\xi
-\frac{1}{2} \int\
\xi\wedge\hat K'(\hat I\xi)\wedge\xi, \label{eq:prolong}
\end{equation}
{\it where the two differential operators $\hat I$, $\hat K$ are the
coefficients of the bivectors in their canonical form.}
\medskip

{\it Proof.}
Let us consider the Schouten-Nijenhuis bracket for the two bivectors and
without loss of generality
take them in the canonical form
\[
\Lambda=\frac{1}{2} \int\theta^{(L)}\xi_A\wedge I^{\langle L\rangle N}_{AB}
D_N\xi_B,
\]
\[
\Psi=\frac{1}{2} \int\theta^{(M)}\xi_C\wedge K^{\langle M\rangle P}_{CD}
D_P\xi_D,
\]
where $\xi_A={\delta}/{\delta\phi_A}$ and operators $\hat I$ , $\hat K$
are antisymmetric. Then we have
\[
\d\Lambda=\frac{1}{2} \int\theta^{(L)}\frac{\partial
I^{\langle L\rangle N}_{AB}}
{\partial\phi_E^{(J)}}\delta\phi_E^{(J)}\xi_A\wedge D_N\xi_B
\]
and
\[
\d\Lambda\inprod\Psi=\frac{1}{4} \int\theta^{(L+M)}
\frac{\partial I^{\langle L\rangle N}_{AB}}{\partial\phi_C^{(J)}}D_J
\left( K^{\langle M\rangle P}_{CD}D_P\xi_D\right)\wedge
\xi_A\wedge D_N\xi_B-
\]
\[
-\frac{1}{4} \int\theta^{(L+M)}D_P\left(
\frac{\partial I^{\langle L\rangle N}_{AB}}
{\partial\phi_D^{(J)}}\xi_A\wedge D_N\xi_B\right)
\wedge D_J\left(\xi_C K^{\langle M\rangle P}_{CD}\right).
\]
Now let us make integration by parts in the second term
\[
\d\Lambda\inprod\Psi=-\frac{1}{4} \int\theta^{(L+M)}\xi_A\wedge
\left(I^{\langle L\rangle N}_{AB}\right)'\left(\hat K^{\langle M\rangle}
\xi\right)\wedge  D_N\xi_B-
\]
\[
-\frac{1}{4} \int\theta^{(L+M+Q)}(-1)^{|P|}{P\choose Q}
\frac{\partial I^{\langle L\rangle N}_{AB}}
{\partial\phi_D^{(J)}}\xi_A\wedge D_N\xi_B
\wedge D_{J+P-Q}\left(\xi_C K^{\langle M\rangle P}_{CD}\right).
\]
At last we change the order of multipliers under wedge product
in the second term,
make a replacement $M\rightarrow M-Q$ and organize the
whole expression in the form
\[
\d\Lambda\inprod\Psi=-\frac{1}{4} \int\theta^{(L+M)}\xi_A\wedge
\left(I^{\langle L\rangle N}_{AB}\right)'_C\Biggl(\hat K^{\langle
M\rangle}_{CD}\xi_D+
\]
\[
+(-1)^{|P|}{P\choose Q}{P-Q\choose R}
D_{P-Q-R}K^{\langle M-Q\rangle P}_{CD}D_R\xi_C \Biggr)\wedge D_N\xi_B.
\]
Having in mind the definition of adjoint operator (\ref{eq:adj}) we can
represent the final result of the calculation as follows,
\[
\bigl[ \Lambda,\Psi\bigr]_{SN}=-\frac{1}{2} \int\theta^{(L+M)}
\xi\wedge\left((\hat I^{\langle L\rangle})'(\hat K^{\langle M\rangle}\xi)
+(\hat K^{\langle M\rangle})'(\hat I^{\langle L\rangle}\xi)\right)\wedge\xi,
\]
thus supporting in this extended formulation the method, proposed
in Ref.~\cite{Olv} for testing the Jacobi identity (see Section 7).
For the general procedure of testing Hamiltonian properties see also
\cite{AstVin}.

\section{Poisson brackets and Hamiltonian vector fields}
Let us call the bivector
\[
\Psi=\frac{1}{2} \int\frac{\delta}{\delta\phi_A}\wedge\hat I_{AB}
\frac{\delta}{\delta\phi_B},
\]
formed with the help of the graded antisymmetric differential operator
\[
\hat I_{AB}=  \theta^{(L)}I^{\langle L\rangle N}_{AB}D_N,
\]
the {\it Poisson bivector} if
\[
\bigl[ \Psi,\Psi\bigr]_{SN} =0.
\]
The operator $\hat I_{AB}$ is then called the {\it Hamiltonian operator}.
We call the value of the Poisson bivector on the differentials of
two functionals $F$ and  $G$
\[
\{ F,G\} = \Psi (\d F,\d G)=\d G\inprod\d F\inprod\Psi
\]
the {\it Poisson bracket} of these functionals.

The explicit form of the Poisson bracket can easily be obtained. It
depends on the explicit form of the functional differential,
which can be changed by the formal partial integration. Of course, all the
possible forms are equivalent. Taking the extreme cases we get an
expression through Fr\'echet derivatives
\begin{equation}
\{ F,G \} =
 \int\theta^{(J)} {\Tr}\biggl( f'_A\hat I^{\langle J\rangle}_{AB}
g'_B \biggr) \label{eq:brack1}
\end{equation}
or through higher Eulerian operators (\ref{eq:higher})
\begin{equation}
\{ F,G \} =
 \int\theta^{(J)} D_{P+Q}\left( E^P_A(f)\hat I^{\langle J\rangle}_{AB}
E^Q_B(g)
\right).\label{eq:brack2}
\end{equation}

\medskip
{\bf Theorem 6.1}
{\it The Poisson bracket defined above satisfies the standard requirements
of the bilinearity, antisymmetry and closeness on the space of local
functionals} {\cal A}, {\it i.e.,}
Definition 2.3 {\it of} \cite{Sol1}.
\medskip

{\it Proof.}
1) From the previous formulas (\ref{eq:brack1}), (\ref{eq:brack2})
it is clear that $\{ F,G \}$ is a local functional,
2) antisymmetry of $\{ F,G \}$ is evident and
3) equivalence of the Jacobi identity
to the Poisson bivector property will be proved in Section 7.
\medskip

The result of interior product of the differential of a local
functional $H$ on the Poisson bivector (up to the sign)
will be called the {\it  Hamiltonian
vector field} (or the {\it Hamiltonian 1-vector})
\[
\hat I\d H=-\d H\inprod\Psi
\]
corresponding to the Hamiltonian $H$.
Evidently, the standard relations take place
\[
\{ F,H\} = \d F(\hat I \d H)=(\hat I \d H)F.
\]

\medskip
{\bf Theorem 6.2}
{\it The Hamiltonian vector field corresponding to the Poisson bracket
of the functionals $F$ and $H$ coincides up to the sign with
commutator of the Hamiltonian vector fields corresponding to these
functionals.}
\medskip

{\it Proof.}
Consider a value of the commutator of Hamiltonian vector fields
$\hat I\d F$ and $\hat I\d H$ on the arbitrary functional $G$
\[
[\hat I\d F, \hat I\d H]G=\hat I\d F(\hat I\d H(G))-\hat I\d H(\hat I\d F(G))
=
\]
\[
=\hat I\d F(\{G,H\})-\hat I\d H(\{G,F\})=\{\{G,H\},F\}-\{\{G,F\},H\}=
\]
\[
=-\{G,\{F,H\}\}=-\hat I\d\{F,H\}(G),
\]
where we have used the Jacobi identity and antisymmetry of Poisson bracket.
Due to the arbitrariness of $G$ the proof is completed.

\medskip
{\it Example 6.3}

Let us consider the first structure
\[
\{ u(x),u(y)\} =\frac{1}{2}(D_x-D_y)\delta(x,y)
\]
of the Korteweg-de Vries equation (Example 7.6 of Ref.~\cite{Olv})
\[
u_t=u_{xxx} +uu_x.
\]
Construct the adjoint graded operator
to $\theta D$ according to Eq.(\ref{eq:adj})
\[
(\theta D)^{\ast}=-\theta D -D\theta
\]
and the antisymmetric operator is
\[
\hat I=\frac{1}{2}\biggl(\theta D- (\theta D)^{\ast}\biggr)=\theta D +
\frac{1}{2}D\theta.
\]
The Poisson bivector has a form
\[
\Psi=\frac{1}{2}\int\theta\biggl( \frac{\delta}{\delta u}\wedge
D\frac{\delta}{\delta u}\biggr).
\]
The differential of a local functional $H$ (for simplicity
it is written in canonical
\[
H=\int\theta h
\]
form) is equal to
\[
\d H=\int\theta h'(\delta u)=
\int\theta^{(k)}(-1)^kE^k(h)\delta u,
\]
where the Fr\'echet derivative or higher Eulerian operators can be used.
Therefore, the Hamiltonian vector field generated by $H$ is
\[
\hat I\d H=-\d H\inprod\Psi=
-\frac{1}{2}\int\theta\left( h'\bigl( D\frac{\delta}{\delta u}
\bigr) - Dh'\bigl( \frac{\delta}{\delta u}\bigr)\right),
\]
or
\[
-\frac{1}{2}\int\theta^{(k)}(-1)^k\left( E^k(h)D-
DE^k(h)\right)\frac{\delta}{\delta u},
\]
or also
\[
-\frac{1}{2}\int\theta^{(k)}(-1)^kD_i\left( E^k(h)D-
DE^k(h)\right)\frac{\partial}{\partial u^{(i)}}.
\]
The value of this vector field on another functional $F$ coincides with
the Poisson bracket
\[
-\d F\inprod\d H\inprod\Psi=\{ F,H\}=
\frac{1}{2} \int\theta^{(k+l)}(-1)^{k+l}
\left( E^k(f)DE^l(h)-E^k(h)DE^l(f)\right).
\]

\section{Proof of Jacobi identity}
In this section we will prove that the Jacobi identity for the Poisson bracket
is fulfilled if and only if the Schouten--Nijenhuis bracket of the
corresponding Poisson bivector with itself is equal to
zero. This should complete the proof of Theorem 6.1.

Let us use one of the possible forms of the Poisson brackets given in
Appendix of \cite{Sol1}
\[
\{ F,G\} =\frac{1}{2} \int\theta^{(J)}{\Tr}\left( f'(\hat I^{\langle
J\rangle}g')
-g'(\hat I^{\langle J\rangle}f')\right),
\]
where the differential operator $\hat I$ is not supposed to be antisymmetric
for the easier comparison of this proof with that given in \cite{Sol1}.
We remind  that in less condensed notations
\[
{\Tr}\left( f'(\hat Ig')\right)= {J\choose M}{K\choose L}
D_L\frac{\partial f}{\partial\phi_A^{(J)}}D_{J+K-L-M}I^N_{AB}
D_{N+M}\frac{\partial g}{\partial\phi_B^{(K)}}
\]
(in Appendix of \cite{Sol1} the indices $M$ and $L$ in the binomial coefficients
of the same formula  are unfortunately given in the opposite order).

We will estimate the bracket
\[
\{\{ F,G\} ,H\} =\frac{1}{2} \int\theta^{(J)}{\Tr}\left(
{\{ f,g\} }'(\hat I^{\langle J\rangle} h')-h'(\hat I^{\langle J\rangle}
{\{ f,g\}}')\right),
\]
where $\{ f,g\}$ denotes the integrand of $\{ F,G\}$. Since Fr\'echet
derivative is a derivation we have
\[
{\{ f,g\} }'=\frac{1}{2} \theta^{(K)}{\Tr}\biggl(
f''(\hat I^{\langle K\rangle}g',\cdot)
+f'\hat I'^{\langle K\rangle}(\cdot)g'+
g''(f'\hat I^{\langle K\rangle},\cdot)-(f\leftrightarrow g)
\biggr)
\]
and
\[
{\Tr}\left( {\{ f,g\}}'\hat Ih'\right)=\frac{1}{2}\left(
f''(\hat Ig',\hat Ih')+f'\hat I'(\hat Ih')g'+g''(f'\hat I,\hat Ih')-
(f\leftrightarrow g)\right).
\]
Let us explain that $f''$ denotes the second Fr\'echet derivative, i.e.,
the symmetric bilinear operator arising in calculation of the second
variation of the local functional $F$ (in the canonical form):

\[
f''(\xi,\eta)= \frac{\partial^2f}
{\partial\phi_A^{(J)}\partial\phi_B^{(K)}}D_J\xi_AD_K\eta_B.
\]
When we  put into entries of $f''$ operators under the trace sign
it should be understood that these operators act on everything except
their own coefficients, for example,
\[
{\Tr}\left(f''(\hat Ig',\hat Ih')\right)
= {L\choose P}{L-P\choose Q}{M\choose S}
{M-S\choose T}\times
\]
\[
\times D_{L+M-P-Q-S-T}\frac{\partial^2f}{\partial\phi_A^{(J)}
\partial\phi_B^{(K)}}
D_{J+T}\left( D_P\hat I_{AC}\frac{\partial g}
{\partial\phi_C^{(L)}}\right) D_{K+Q}\left( D_S\hat I_{BD}\frac{\partial h}
{\partial\phi_D^{(M)}}\right)
\]
and the expression remains symmetric under permutation of its entries
\[
{\Tr}\left( f''(\hat Ig',\hat Ih')\right)=
{\Tr}\left( f''(\hat Ih',\hat Ig')\right) .
\]
When the operator $\hat I$  stands to the right from the operator
of Fr\'echet derivative $f'$ as in expression
\[
{\Tr}\left( g''(\hat Ih',f'\hat I)\right) ,
\]
it acts on everything except $f'$. At last, for Fr\'echet derivative
of the operator we have
\[
\hat I'(\hat Ih')= \frac{\partial I^K_{AB}}{\partial\phi_C^{(J)}}
D_J\left( I^L_{CD}D_L\frac{\partial h}{\partial\phi_D^{(M)}}D_M\right)D_K.
\]
Making similar calculations we get
\[
{\Tr}\left( h'\hat I{\{ f,g\}}'\right)=
\frac{1}{2}{\Tr}\left( f''(h'\hat I,\hat Ig')+f'\hat I'(h'\hat I)g'+
g''(f'\hat I,h'\hat I)-(f\leftrightarrow g)\right)
\]
and therefore
\[
\{\{ F,G\} ,H\} =\frac{1}{4} \int\theta^{(J+K)}{\Tr}\Biggl(
f''(\hat I^{\langle J\rangle}g',\hat I^{\langle K\rangle}h')-
f''(h'\hat I^{\langle J\rangle},\hat I^{\langle K\rangle}g')-
\]
\[
-f''(\hat I^{\langle J\rangle}h',g'\hat I^{\langle K\rangle})+
f''(g'\hat I^{\langle J\rangle},h'\hat I^{\langle K\rangle})+
f'\hat I'^{\langle J\rangle}(\hat I^{\langle K\rangle}h'-h'
\hat I^{\langle K\rangle})g' -(f\leftrightarrow g)
\Biggr) .
\]
Just the first four terms, apart from the fifth containing Fr\'echet
derivative of the
operator $\hat I$, were present in our proof for nonultralocal case
given in \cite{Sol1} (only terms with zero grading were allowed for $\hat I$ there).
After cyclic permutation of $F$, $G$, $H$ all terms with the symmetric
operator of the second Fr\'echet derivative are mutually cancelled and
\[
\{\{ F,G\} ,H\} + {\rm c.p.}=\frac{1}{4}\int\theta^{(J+K)}{\Tr}\Biggl(
f'\hat I'^{\langle J\rangle}(\hat I^{\langle K\rangle}h'-
h'\hat I^{\langle K\rangle })g'-
\]
\[
-g'\hat I'^{\langle J\rangle}
(\hat I^{\langle K\rangle}h'-h'\hat I^{\langle K\rangle})f'+ c.p.
\Biggr),
\]
where cyclic permutations of $F$, $G$, and $H$ are abbreviated to $c.p.$.
When operator $\hat I$ is given in explicitly antisymmetric form all
the four terms are equal. Taking into account Olver's Lemma
(\ref{eq:prolong})
we get
\[
\{\{ F,G\} ,H\} +  c.p. =-\bigl[\hat I,\hat
I\bigr]_{SN}(\d F,\d G,\d H),
\]
so finishing the proof.

\section{Examples of nonultralocal operators}
The second structure of the Korteweg-de Vries equation may serve as a
counter-example to the hypothesis \cite{coll}
that all operators which are Hamiltonian
with respect to the standard Poisson brackets should also be Hamiltonian in
the new brackets.

\medskip
{\it Example 8.1}

Let us start with the standard expression (Example 7.6 of Ref.~\cite{Olv})
\[
\{ u(x),u(y)\} =\left(\frac{d^3}{dx^3}+\frac{2}{3}u\frac{d}{dx}+
\frac{1}{3}\frac{du}{dx}\right)\delta(x,y)
\]
and construct the adjoint operator to
\[
\hat K=\theta\left(D_3+\frac{2}{3}uD+\frac{1}{3}Du\right),
\]
which is
\[
\hat K^{\ast}=-\theta\left(D_3+\frac{2}{3}uD+\frac{1}{3}Du\right)
-D\theta\left(3D_2+
\frac{2}{3}u\right)-3D_2\theta D-D_3\theta.
\]
Then the antisymmetric operator
\[
\hat I=\frac{1}{2}(\hat K-\hat K^{\ast})=
\theta
\left(D_3+\frac{2}{3}uD+\frac{1}{3}Du\right)
+D\theta\left(\frac{3}{2}D_2+\frac{1}{3}u\right)
+\frac{3}{2}D_2\theta D+\frac{1}{2}D_3\theta
\]
can be used for forming the bivector
\[
\Psi=\frac{1}{2}\int\xi\wedge\hat I\xi,
\]
where  ${\delta}/{\delta u}=\xi$. This bivector has a form
\[
\Psi=\frac{1}{2}\int\left(\theta\xi\wedge D_3\xi+\frac{3}{2}D\theta\xi
\wedge D_2\xi+(\frac{3}{2}D_2\theta+\frac{2}{3}\theta u)\xi\wedge D\xi\right).
\]
Then evaluating the Schouten-Nijenhuis bracket for the bivector with
the help of Statement 5.3
\[
\bigl[\Psi,\Psi\bigr]_{SN}=\int\left(\frac{2}{3}\theta\xi\wedge D_3\xi
\wedge D\xi +
D\theta\xi\wedge D_2\xi\wedge D\xi\right)
\]
and  integrating the first term by parts we get
\[
\bigl[\Psi,\Psi\bigr]_{SN}=\frac{1}{3}\int\theta D\bigl(\xi\wedge D\xi
\wedge D_2\xi\bigr).
\]
Therefore, instead of the Jacobi identity we  have
\[
\{\{ F,G\},H\}+c.p.=-\frac{1}{3}
\int\limits_{\Omega}D_{i+j+k+1}\left( E^i(f)DE^j(g)D_2E^k(h)+c.p.\right)dx.
\]
So, the second structure of KdV equation can be Hamiltonian only
under special boundary conditions.

\medskip
{\it Example 8.2}

Now let consider another example which is also nonultralocal, but the operator
remains
to be Hamiltonian in the new brackets independently of boundary conditions.
The Euler equations for the flow of ideal  fluid
can be written \cite{Olv} in Hamiltonian form as follows
(Example 7.10 of Ref.~\cite{Olv})
\[
\frac{\partial{\bf\omega}}{\partial t}={\cal D}\frac{\delta H}
{\delta{\bf\omega}},
\]
where
\[
H=\int\frac{1}{2}\vert {\bf u}\vert^2d^2x,\qquad {\bf\omega}={\bf\nabla}
\times{\bf u}.
\]
Let us limit our consideration by the 2-dimensional case when $\bf\omega$
has only one component $\omega$ and
\[
{\cal D}={\bf\omega}_xD_y-{\bf\omega}_yD_x,
\]
where $\omega_i=D_i\omega$, $i=(x,y)$.
We can construct the antisymmetric operator
\[
\hat I=\frac{1}{2}\biggl( \theta{\cal D}-(\theta{\cal D})^{\ast}\biggr)=
\theta(\omega_xD_y-\omega_yD_x)
+\frac{1}{2}(D_y\theta\omega_x-D_x\theta\omega
_y),
\]
and then the bivector
\[
\Psi=\frac{1}{2}\int\xi\wedge\hat I\xi=
\frac{1}{2}\int\theta(\omega_x\xi\wedge
\xi_y-\omega_y\xi\wedge\xi_x),
\]
where $\xi={\delta}/{\delta\omega}$.
Statement 5.3 gives us
\[
\bigl[\Psi,\Psi\bigr]_{SN}=
\int\Biggl(\theta\left(
\omega_x(\xi\wedge\xi_{xy}\wedge\xi_y-\xi\wedge\xi_{yy}
\wedge\xi_x)+\omega_y(\xi\wedge\xi_{xy}\wedge\xi_x-\xi\wedge\xi_{xx}\wedge
\xi_y)\right)+
\]
\[
+\left(D_y\theta\omega_x-
D_x\theta\omega_y\right)\xi\wedge\xi_x\wedge\xi_y\Biggr)
\]
and after integration by parts the expression can be reduced to zero.

\section{Conclusion}
We have shown that there is an extension of the standard formal
variational calculus which incorporates  divergences without
any specification of  boundary conditions. 
It should be important to understand relations of this
formalism to the constructions of the variational bicomplex \cite{And}.
It seems also rather interesting to study if some physically relevant
algebras can be realized with the help of the new Poisson brackets
as algebras of local functionals. One such example is considered
in \cite{Sol97}.

\section*{Acknowledgements}

It is a pleasure to thank I.~Kanatchikov, J.~Kijowski, O.~Mokhov, J.~Nester,
C.~Rovelli, S.~Storchak and E.~Timoshenko for discussions on this work.

\end{document}